\input amstex
\documentstyle{amsppt}
%%%%%%%%%%%%%%%%%%%%%%%%
%\input coverdef
 \loadmsbm
  \loadbold

 \magnification=\magstep1
        \pagewidth{13cm}
        \pageheight{20cm}
\topmatter

\define\dN{{\Bbb N}}

\define\dQ{{\Bbb Q}}

\define\dZ{{\Bbb Z}}

\define\cF{{\Cal F}}

\define\ep{{\frak p}}
\define\eP{{\frak P}}

\define\tw#1{\,^{#1}\!}

\font\pfeile = cmsy10 scaled 1440
\newfam\pfeilfam
\textfont\pfeilfam=\pfeile
                   \scriptfont\pfeilfam=\pfeile
                                      \scriptscriptfont\pfeilfam=\pfeile

\mathchardef\swpfeil="2D2E

\def\\{\let\stoken= } \\
\long\def\unexpandedwrite#1#2{\def\finwrite{\write#1}%
{\aftergroup\finwrite\aftergroup{\sanitize#2\endsanity}}}
\def\sanitize{\futurelet\next\sanswitch}
\def\sanswitch{\ifx\next\endsanity
\else\ifcat\noexpand\next\stoken\aftergroup\space\let\next=\eat
\else\ifcat\noexpand\next\bgroup\aftergroup{\let\next=\eat
\else\ifcat\noexpand\next\egroup\aftergroup}\let\next=\eat
\else\let\next=\copytoken\fi\fi\fi\fi \next}
\def\eat{\afterassignment\sanitize \let\next= }
\long\def\copytoken#1{\ifcat\noexpand#1\relax\aftergroup\noexpand
\else\ifcat\noexpand#1\noexpand~\aftergroup\noexpand\fi\fi
\aftergroup#1\sanitize}
\def\endsanity\endsanity{}

\expandafter\ifx\csname pre numero.tex at\endcsname\relax \else
\endinput\fi \expandafter\chardef\csname pre numero.tex
at\endcsname=\the\catcode`\@ \catcode`\@=11
%%%%%%%%%%%
\def\tokenitize{\futurelet\next\tokenswitch}
\def\tokenswitch{\ifx\next\endtokenity
\else\ifcat\noexpand\next\stoken\aftergroup\space\let\next=\eatt
\else\let\next=\copytok\fi\fi\next}
\def\eatt{\afterassignment\tokenitize \let\next= }
\long\def\copytok#1{\aftergroup\string\aftergroup#1\tokenitize}
\def\endtokenity{}
\newcount \secno
\newcount \Refno
\def\phantomlabel#1{{\aftergroup\expandafter\aftergroup\xdef%
\aftergroup\csname\tokenitize ref#1\endtokenity\aftergroup\endcsname}%
{{\the\secno.\the\Refno}}%
\global\advance\Refno by 1\relax}
\def\label#1{{\the\secno.\the\Refno}\phantomlabel{#1}%
\unexpandedwrite\labx{\phantomlabel{#1}}}
\def\Ref#1{{\aftergroup\expandafter\aftergroup\ifx
\aftergroup\csname\tokenitize ref#1\endtokenity}\endcsname\relax
??R\'ef\'erence {\tokenitize#1\endtokenity} non d\'efinie??%
\message{Reference #1 non definie}%
\else{\aftergroup\csname\tokenitize
ref#1\endtokenity}\endcsname\fi
\ifnextsecn@\else\message{ATTENTION secno n'a pas ete
initialise}\fi }
\def\phantomsection{\global\advance\secno by 1 \global\Refno = 1\relax}
\def\secnum{{\phantomsection\the\secno\unexpandedwrite\labx{\phantomsection}}}
\newif\ifnextsecn@
\def\nextsecno#1{\nextsecn@true\global\secno=#1\global\advance\secno by -1%
\immediate\write\labx{\secno=\the\secno}\message{\the\secno} }
\openin1 \jobname.lab \ifeof1 \message
                                {No file \jobname.lab.}\else\closein1\relax\input \jobname.lab \fi
\newwrite\labx
\immediate\openout\labx=\jobname.lab %nom du fichier de labels
%%%%%%%%%%%%
%  Restore the catcode value for @ that was previously saved.
\catcode`\@=\csname pre numero.tex at\endcsname
%%%%%%%%%%%%
\nextsecno 1 \scrollmode \NoBlackBoxes
\define\BeBi{1}
\define\Br{2}
\define\BBH{3}
\define\Hu{4}
\define\Me{5}
\define\MM{6}
\define\NSW{7}
\define\So{8}
\define\Soa{9}
\title
Polynomials with roots in $\dQ_p$ for all $p$
\endtitle
\author
 Jack Sonn \endauthor
\affil  Technion--Israel Institute of Technology, Haifa, Israel
\endaffil

\address
Department of Mathematics, Technion, 32000 Haifa, Israel
\endaddress \email
 sonn\@math.technion.ac.il
\endemail

%\thanks

%The research  was supported by
%the Fund for the Promotion of
%Research at the Technion.

%\endthanks
\abstract Let $f(x)$ be a monic polynomial in $\dZ[x]$ with no
rational roots but with roots in $\dQ_p$ for all $p$, or
equivalently, with roots mod $n$ for all $n$. It is known that
$f(x)$ cannot be irreducible but can be a product of two or more
irreducible polynomials, and that if $f(x)$ is a product of $m>1$
irreducible polynomials, then its Galois group must be a union of
conjugates of $m$ proper subgroups.  We prove that for any $m>1$,
every finite solvable group which is  a union of conjugates of $m$
proper subgroups (where all these conjugates have trivial
intersection) occurs as the Galois group of such a polynomial, and
that the same result (with $m=2$) holds for all Frobenius groups.
It is also observed that every nonsolvable Frobenius group is
realizable as the Galois group of a geometric--i.e. regular--
extension of $\dQ(t)$.
\endabstract \endtopmatter

\head
%%%%%
  \secnum.  Introduction
\endhead
\smallskip

\document
Let $f(x)$ be a monic polynomial with rational integer
coefficients.  Suppose $f$ has no rational roots but has a root in
$\dQ_p$ for all $p$, or equivalently, has a root mod $n$ for all
positive integers $n$. It has been observed \cite\BeBi,\cite\Br
that $f$ cannot be irreducible (even under the slightly weaker
assumption that $f$ has a root mod $p$ for all prime numbers $p$),
and that if $f$ is a product of $m>1$ irreducible polynomials,
then its Galois group must be a union of conjugates of $m$ proper
subgroups. \footnote {In \cite\BeBi,  an effective condition is
given for checking if $f$ has this property. }
 Indeed, if $f$ is irreducible,
let $\alpha$ be a root of $f$, $F=\dQ(\alpha)$, $K$ the splitting
field of $f$ over $\dQ$, $G$ the Galois group $G(K/\dQ)$.  Let $p$
be a prime not dividing the discriminant $\text{disc}(f)$ of $f$.
$f$ has a root mod $p$ implies that for some prime $\ep$ of $K$
dividing $p$, the decomposition group $G(\ep)$ of $\ep$ is
contained in $G(K/F)$, or equivalently, for any prime $\ep$ of $K$
dividing $p$, there exists a root $\beta$ of $f$ such that
$G(\ep)\subseteq G(K/\dQ(\beta))$.  By Chebotarev's density
theorem, $G=\bigcup_{\ep\nmid\text{disc}(f)}G(\ep)$.  But then
$$\bigcup_{\ep\nmid\text{disc}(f)}G(\ep)\subseteq
\bigcup_{f(\beta)=0}G(K/\dQ(\beta))$$ cannot equal $G$ since a
union of conjugates of a proper subgroup of a finite group cannot
equal the whole group.  Thus $f$ cannot be irreducible.  By a
similar argument, if  $f$ is a product $g_1\cdots g_m$ of $m$
irreducible polynomials in $\dZ[x]$ (of degree bigger than $1$)
with splitting field $K$ and Galois group $G$, then $G$ is a union
of conjugates of $m$ proper subgroups, namely
$G(K/\dQ(\alpha_i))$, where $g_i(\alpha_i)=0$, $i=1,...,m$.
Finite groups which are the union of conjugates of two proper
subgroups are studied in \cite \Br, where it is proved among other
things that if a group $G$ has this property, and the two
subgroups are nilpotent, then $G$ is solvable.  It is also proved
in \cite \BBH that if the symmetric group $S_n$ has this property
then $3\leq n\leq 6$.  In \cite \Br (see also \cite \BeBi), the
polynomials $(x^r-2)\Phi_r(x)$, $r\geq 3$ a prime ($\Phi_r(x)$ is
the $r$th cyclotomic polynomial), with Galois groups the Frobenius
groups of order $r(r-1)$, are given as examples of polynomials
with no rational roots and roots mod $p$ for all $p$, and the
inverse problem is raised: if $G$ is a union of conjugates of two
proper subgroups (one should also assume that the intersection of
all these conjugates is trivial), then can $G$ be realized as the
Galois group of the product $f$ of two irreducible nonlinear
polynomials such that $f$ has a root mod $p$ for all $p$?  The
answer appears not to have been known even for the dihedral group
of order ten.  We will prove that every finite solvable group with
the above property can be realized in this way, with $f$ having a
root in $\dQ_p$ for all $p$. As it turns out, Shafarevich's
realization of solvable groups already yields extensions with the
required property (even for ``$m$ proper subgroups" and ``$m$
irreducible polynomials" instead of two). We will also prove the
result for all nonsolvable Frobenius groups. On the other hand,
the question seems to be open even for the symmetric group $S_6$,
which as we mentioned above, is the union of conjugates of two
proper subgroups.  In connection with Frobenius groups, we observe
that every nonsolvable Frobenius group is realizable as the Galois
group of a geometric--i.e. regular-- extension of $\dQ(t)$, a fact
that does not seem to have been pointed out before.

\bigskip

\head \secnum.
Roots in $\dQ_p$ for all $p$ \endhead
\medskip
We begin with a characterization of Galois extensions which are
splitting fields of polynomials which are products of $m$
irreducible nonlinear polynomials in $\dQ[x]$ and which have roots
in $\dQ_p$ for all $p$.  Note that if $f \in \dZ[x]$ has a root in
$\dQ_p$ then $f$ has a root mod $p$.

%\proclaim {Lemma \label{basic}}

\proclaim {Proposition \label{char}} Let $K/\dQ$ be a finite
Galois extension with Galois group $G$.  The following are
equivalent:
\smallskip

(1) \ \ \ \ \ $K$ is the splitting field of a product $f=g_1
\cdots g_m$ of $m$ irreducible polynomials of degree greater than
$1$ in $\dQ[x]$ and $f$ has a root in $\dQ_p$ for all primes $p$.
\smallskip

(2) \ \ \ \ \ $G$ is the union of the conjugates of $m$ proper
subgroups $A_1,...,A_m$, the intersection of all these conjugates
is trivial, and for all primes $\ep$ of $K$, the decomposition
group $G(\ep)$ is contained in a conjugate of some $A_i$.
\endproclaim

\bigskip

\demo{Proof}  Assume first that (1) holds, i.e. $K$ is the
splitting field of a product $f=g_1 \cdots g_m$ of $m$ irreducible
polynomials of degree greater than $1$ and $f$ has a root in
$\dQ_p$ for all primes $p$.  Let $\alpha_1,...,\alpha_m$ be roots
of $g_1,...,g_m$ respectively in $K$ and let
$A_i:=G(K/\dQ(\alpha_i)), 1\leq i \leq m$.  Let $p$ be a given
prime number. By assumption $f$ has a root in $\dQ_p$, hence some
$g_i$ has a root in $\dQ_p$.  Then for some prime $\ep$ of $K$
dividing $p$, the decomposition group $G(\ep)$ is contained in
$A_i$, or equivalently, for every prime $\ep$ of $K$ dividing $p$,
the decomposition group $G(\ep)$ is contained in some conjugate
$G(K/\dQ(\alpha_i'))$ of $A_i$.   We therefore conclude that if
$f$ has a root in $\dQ_p$ for all $p$, then for all primes $\ep$
of $K$, the decomposition group $G(\ep)$ is contained in a
conjugate of some $A_i$.  By Chebotarev's density theorem, every
cyclic subgroup of $G$ occurs as a decomposition group of some
(unramified) prime $\ep$ of $K$, hence $G$ is the union of the
conjugates of $A_1,...,A_m$.  The intersection of all the
conjugates of $A_1,...,A_m$ is trivial because $K$ is the
splitting field of $f$.  Thus (2) holds.

Conversely, assume (2), i.e. $G$ is the union of the conjugates of
$m$ proper subgroups $A_1,...,A_m$, the intersection of all these
conjugates is trivial, and for all primes $\ep$ of $K$, the
decomposition group $G(\ep)$ is contained in a conjugate of some
$A_i$.  Let $\alpha_1,...,\alpha_m \in K$ such that
$A_i=G(K/\dQ(\alpha_i)$.  Let $g_1,...,g_m$ be the minimal
polynomials of $\alpha_1,...,\alpha_m$ resp. over $\dQ$.  Since
the intersection of the conjugates of $A_1,...,A_m$ is trivial,
$K$ is the splitting field of $f=g_1 \cdots g_m$ over $\dQ$, and
since $A_1,...,A_m$ are proper subgroups of $G$, $g_1,...,g_m$
have degrees greater than $1$.  Let $p$ be a rational prime, $\ep$
a prime of $K$ dividing $p$.  By assumption the decomposition
group $G(\ep)$ of $\ep$ is contained in a conjugate of some $A_i$.
Then $g_i$ has a root in $\dQ_p$.  Thus our assumptions imply that
$f$ has a root in $\dQ_p$ for all $p$. \qed\enddemo

Note that (2) implies that $G$ is necessarily noncyclic.
\smallskip
 We now prove a realization theorem for solvable groups.

\proclaim {Theorem \label{solvable}}  Let $G$ be a finite solvable
group which is the union of the conjugates of $m$ proper
subgroups, where the intersection of all these conjugates is
trivial. Then there exists a polynomial $f(x)$ which is the
product of $m$ irreducible nonlinear polynomials in $\dQ[x]$ with
Galois group $G$ and having a root in $\dQ_p$ for all rational
primes $p$.  In particular {\rm{(since every noncyclic group is a
union of (conjugates of) proper subgroups with trivial
intersection)}}, every noncyclic finite solvable group is
realizable as the Galois group over $\dQ$ of a polynomial $f(x)
\in \dQ[x]$ having no rational roots and having a root in $\dQ_p$
for all rational primes $p$.
\endproclaim

\demo{Proof}  The proof will follow easily from the observation
that Shafarevich's realization of solvable groups as Galois groups
over number fields yields an extension $K/\dQ$ with all
decomposition groups $G(\ep)$ cyclic.  Indeed, let $G$ be a finite
solvable group which is the union of the conjugates of $m$ proper
subgroups, and suppose $K/\dQ$ is Galois with group $G$ with all
decomposition groups $G(\ep)$ cyclic.  Then every decomposition
group is contained in a conjugate of some $A_i$, so by
Proposition~\Ref{char}, we are done.

To verify the observation about Shafarevich's realization of
solvable groups, we use the exposition of the proof of
Shafarevich's theorem in \cite\NSW. The key result in the
construction is \cite{\NSW, Theorem 9.5.11}. Let $\cF(n)$ denote
the free pro-$p-G$ operator group on $n$ generators. $G$ acts
``freely" on $\cF(n)$.  There is a filtration $\cF(n)^{(\nu)}$
($\nu\in \dN\times\dN$) defined on $\cF(n)$ which is a refinement
of the descending $p$-central series, all of whose terms are
$G$-invariant.  (For the precise definition see \cite{\NSW, pp.
481ff}.)   Now \cite{\NSW, Theorem
 9.5.11} says that if $K/k$ is any Galois extension of global
 fields with group $G$, then for any $p,n,\nu$, the split embedding
 problem associated with the epimorphism $\cF(n)/\cF(n)^{(\nu)}\rtimes G\twoheadrightarrow
 G$ has a proper solution with solution field $N_{\nu}^n$, such
 that  (if $p\neq\text{char}(K)$), all divisors of $p$, all infinite primes,
 and all primes of $K$ which are ramified in $K/k$ split completely in
 $N_{\nu}^n/K$, and all primes $\ep$ of $K$ which ramify in
 $N_{\nu}^n/K$ split completely in $K/k$, and the local extension $N_{\nu,\ep}^n/k_{\ep}$ is a
 totally ramified--hence cyclic--extension.  In particular, if all
 decomposition groups in $G(K/k)$ are cyclic, then all decomposition groups in $N_{\nu}^n/k$ are cyclic.
 Now given any semidirect product $P\rtimes G$, with $P$ a finite
$p$-group, there exists some $n,\nu$
 and an operator epimorphism from $\cF(n)/\cF(n)^{(\nu)}$ to $P$, and a
 corresponding epimorphism from the semidirect product
 $\cF(n)/\cF(n)^{(\nu)}\rtimes G$ to $P\rtimes G$.  This implies
 that the split embedding problem associated with the epimorphism
 $P\rtimes G\twoheadrightarrow G$ has a proper solution with
 solution field $N\subseteq N_{\nu}^n$, and if  all
 decomposition groups in $G(K/k)$ are cyclic, then all decomposition groups in $G(N/k)$ are
 cyclic.  Finally, given any semidirect product $Q\rtimes G$ with
 $Q$ a finite nilpotent group, the above argument implies that the
 embedding problem associated with the epimorphism $Q\rtimes G\twoheadrightarrow G$ has a proper solution with
 solution field $M$, and if  all
 decomposition groups in $G(K/k)$ are cyclic, then all decomposition groups in $M/k$ are
 cyclic.  \smallskip Now the proof of Shafarevich's theorem
 follows by applying a theorem of Ore: let $G$ be a finite
 solvable group.  Then $G$ has a nilpotent normal subgroup $Q$ and
 a proper subgroup $H$ such that $G=QH$.  By induction we may
 assume $H$ is realized as a Galois group $G(K/k)$ with all
 decomposition groups cyclic.  Consider the semidirect product $Q\rtimes H$ with $H$ acting on $Q$ by conjugation
 inside $G$.  By the above, the embedding problem
 associated with the epimorphism $Q\rtimes H\twoheadrightarrow H$ has a proper solution with
 solution field $M$, and since  all
 decomposition groups in $G(K/k)$ are cyclic,  all decomposition groups in $M/k$ are
 cyclic.  Finally, since $G$ is a homomorphic image of $Q\rtimes
 H$, there is a subfield $L$ of $M$  such that $L/k$
 is Galois with group $G$ and all decomposition groups cyclic.
 This verifies the observation about Shafarevich's construction
 and completes the proof of the Theorem. \qed \enddemo
\medskip

\bf Remark. \rm Proposition~\Ref{char} and Theorem~\Ref{solvable}
hold with the base field $\dQ$ replaced by an arbitrary global
field $k$, with the same proof, where the primes $p$ are replaced
by the primes of $k$. \bigskip

We now turn to nonsolvable groups.  One family of groups each of
which is the union of two conjugacy classes of proper subgroups
 is the family of Frobenius
groups.  Unlike most nonsolvable groups, nonsolvable Frobenius
groups are known to be realizable as Galois groups over $\dQ$
\cite\So. If $G$ is a Frobenius group, then $G$ is a semidirect
product $Q\rtimes H$, where $g^{-1}Hg\cap H =1$ for all $g\notin
H$, which implies that $G$ is covered by $Q$ (which is normal) and
the conjugates of $H$ \cite{\Hu, p.495}.

\proclaim{Theorem \label{frob}} Let $G$ be a Frobenius group.
Then there exists a polynomial $f(x)$ which is the product of two
irreducible polynomials in $\dQ[x]$ with Galois group $G$ and
having a root in $\dQ_p$ for all rational primes $p$. \endproclaim

\demo{Proof}
%We may assume $G$ nonsolvable in view of
%Theorem~\Ref{solvable}.
By Proposition~\Ref{char}, it suffices to show that $G$ is
realizable as the Galois group of an extension $K/\dQ$ with each
decomposition group contained in either $Q$ or a conjugate of $H$.
One of the facts about Frobenius groups is very relevant here,
namely Thompson's theorem that the ``Frobenius kernel" $Q$ of $G$
is nilpotent \cite{\Hu, p. 499, Thm. 8.7} (in fact for $G$
nonsolvable, $Q$ is even abelian \cite{\Hu, p. 506, Thm. 8.18}).
We may therefore use Shafarevich's theorem as we did in the proof
of Theorem~\Ref{solvable}.  The same argument that we used there
shows that if we can realize the ``Frobenius complement" $H$ by a
Galois extension $L/\dQ$, then we can embed $L/\dQ$ into a Galois
extension $K/\dQ$ with group $G$, such that the ramified primes in
$L/\dQ$ split completely in $K/L$, and the ramified primes of
$K/L$ are split completely in $L/\dQ$ with cyclic decomposition
groups in $K/\dQ$. Let $\ep$ be a prime of $K$. If it is
unramified over $\dQ$, its decomposition group is cyclic, hence
contained in either $Q$ or a conjugate of $H$.  If it is ramified
over $\dQ$, let $I(\ep)$ be its inertia group. If $I(\ep)\subseteq
Q$, then $\ep$ is unramified in $L/\dQ$ and ramified in $K/L$,
hence split completely in $L/\dQ$, so its decomposition group is
contained in $Q$ (and in fact equals $I(\ep)$) and we are done.
Otherwise, $I(\ep)$ is not contained in $Q$, so $\ep$ is ramified
in $L/\dQ$, and so splits completely in $K/L$.  This means that
its decomposition group $G(\ep)$ intersects $Q$ trivially.  We
will use the following group-theoretic lemma, for which we are
indebted to David Chillag.

\proclaim{Lemma \label{lemma}} Let $G$ be a Frobenius group
$Q\rtimes H$.  Then every subgroup $D$ of $G$ such that $D\cap
Q=\{1\}$ is contained in a conjugate of $H$.
\endproclaim

\demo{Proof}  Since $D\cap Q=\{1\}$, every element of $D$ acts
without fixed points on $Q$, so $DQ$ is a Frobenius group with
kernel $Q$ and complement $D$.  By \cite{\Hu, Thm. 8.18, p. 506},
the center of any Frobenius complement and in particular, $Z(D)$,
is nontrivial.  Let $d$ be a nontrivial element of $Z(D)$.  Then
$D\subseteq C(d)$, where $C(d)$ is the centralizer of $d$ in $G$.
$d$ lies in a conjugate of $H$, so without loss of generality, we
may assume $d\in H$.  On the other hand, $C(d)\subseteq H$, since
if $x\in C(d)\setminus H$, then $xdx^{-1}=d$ lies in $xHx^{-1}\cap
H=\{1\}$, contradiction.  We therefore have $D \subseteq H$.
 \qed \enddemo

By the argument preceding the Lemma, we have

 \proclaim{Corollary  \label {corollary}} Let $G$ be a Frobenius
 group $Q\rtimes H$ and
assume $L/\dQ$ is Galois with $G(L/\dQ)\cong H$.   Then $L/\dQ$
can be embedded into a Galois extension $K/\dQ$ with
$G(K/\dQ)\cong G$, such that for all primes $\eP$ of $K$, the
decomposition group of $\eP$ is contained in $Q$ or a conjugate of
$H$. \endproclaim

The proof of Theorem~\Ref{frob} is then completed by  the
realization of nonsolvable Frobenius complements over $\dQ$ in
\cite\So, or by using the following remark.    \qed \enddemo

\proclaim{ Remark \label{remark}}  Theorem \Ref{frob} holds with
$\dQ$ replaced by an arbitrary number field $k$.
\endproclaim \smallskip

\demo{Proof} It suffices to realize every nonsolvable Frobenius
complement $H$ as the Galois group of a geometric (regular over
$\dQ$) extension of the rational function field $\dQ(t)$, since
Hilbert's Irreducibility Theorem implies that if a group $G$ is
realizable as the Galois group of a geometric extension of
$\dQ(t)$, then it is realizable over every number field. Now $H$
is itself a semidirect product $Z \rtimes B$, where $Z$ is the
semidirect product of two cyclic groups of orders relatively prime
to each other and to $2,3,5$, and $B$ is one of two groups: $\hat
A_5$, the double cover  of the alternating group $A_5$, or $\hat
S_5$, one of the two double covers of the symmetric group $S_5$.
Over any Hilbertian field $F$, every split geometric embedding
problem with abelian kernel has a proper geometric solution
\cite{\MM, Thm. 2.4, p. 275}. Two applications of this fact reduce
the proof to the geometric realization of $\hat A_5$ and $\hat
S_5$ over $\dQ(t)$, which appear in \cite\Me and \cite \Soa,
respectively. \qed \enddemo

 The proof of
Remark~\Ref{remark} implies a result that does not seem to have
been observed before: \proclaim{Theorem \label{geomfrob}} Every
nonsolvable Frobenius group is realizable as the Galois group of a
geometric extension of $\dQ(t)$.
\endproclaim

\demo{Proof} As mentioned earlier, the Frobenius kernel of a
nonsolvable Frobenius group is abelian.  Since the Frobenius
complement is realizable geometrically over $\dQ(t)$, another
application of \cite{\MM, Thm. 2.4, p. 275} yields the result.
\qed \enddemo

Remarkably, this result is not known for solvable Frobenius groups
in general, since nonabelian Frobenius kernels are known to exist.

\enddocument